\documentclass[twocolumn]{autart}    % Enable this line and disable the 
\usepackage{amsmath,amssymb,amsfonts}
\usepackage{graphicx}
\usepackage{textcomp}
\usepackage{mathtools}
\usepackage{import}
\usepackage{xcolor}
\usepackage{tikz}
\usepackage{caption}

\usepackage{savesym}
\savesymbol{AND}
\usepackage{algorithm}
\usepackage[]{algpseudocode}
\usepackage{xcolor}
\usepackage{natbib}        % required for bibliography
\usepackage{booktabs}
\usepackage{enumerate}% http://ctan.org/pkg/enumerate
\usepackage{bbm}

\usepackage{graphicx}          % Include this line if your 
\begin{document}   

\begin{frontmatter}

\title{A stochastic MPC scheme for distributed systems with multiplicative uncertainty
\thanksref{footnoteinfo}} % Title, preferably not more 
                                                % than 10 words.

\thanks[footnoteinfo]{This paper was not presented at any IFAC 
meeting. Corresponding author C.~Mark. \\
© 2022. This manuscript version is made available under the CC-BY-NC-ND 4.0 license https://creativecommons.org/licenses/by-nc-nd/4.0/}

\author[First]{Christoph Mark} \ead{mark@eit.uni-kl.de},
\author[First]{Steven Liu}  \ead{sliu@eit.uni-kl.de}

\address[First]{Institute  of  Control  Systems,  Department  of  Electrical  and  Computer
Engineering, University of Kaiserslautern, Erwin-Schrödinger-Str. 12, 67663
Kaiserslautern, Germany}

                             % keyword list or with the 
                                          % help of the Automatica 
                                          % keyword wizard

\begin{abstract}                          % Abstract of not more than 200 words.
This paper presents a Distributed Stochastic Model Predictive Control algorithm for networks of linear systems with multiplicative uncertainties and local chance constraints on the states and control inputs. The chance constraints are approximated via Cantelli's inequality by means of expected value and covariance. The cooperative control algorithm is based on the distributed Alternating Direction Method of Multipliers, which renders the controller fully distributedly implementable, recursively feasible and ensures point-wise convergence of the states. The aforementioned properties are guaranteed through a properly selected distributed invariant set and distributed terminal constraints for the mean and covariance. The paper closes with an example highlighting the chance constraint satisfaction, numerical properties and scalability of our approach.
\end{abstract}

\end{frontmatter}
\section{Introduction}
Model Predictive Control (MPC) is an optimization based control strategy \cite{kouvaritakis2016model}, which received a lot of attention during the last couple of decades in academic research and industrial applications. One of the main drawbacks of MPC is the computational demand, since in each time step an optimization problem has to be solved. In order to apply MPC to large-scale systems, we have to consider distributed approaches, that is, distributed MPC (DMPC) \cite{christofides2013distributed,conte2016distributed}.
\\
In the presence of stochastic uncertainty we fall in the domain of stochastic MPC \cite{mesbah2016stochastic}. There exist basically two approaches, namely scenario-based methods \cite{blackmore2010probabilistic,bernardini2009scenario} and analytical approximation methods \cite{paulson2020stochastic,hewing2018stochastic,lorenzen2016constraint,
cannon2010stochastic,primbs2009stochastic,cannon2009model}. Scenario-based methods rely on a sufficient number of disturbance realizations in order to compute an optimal solution to the stochastic MPC problem via a sampling-average approximation, whereas in analytical approximation methods the stochastic MPC problem is reformulated as a deterministic one by exploiting the distributional information of the uncertainty. In this work, we consider distributed linear systems with multiplicative noise motivated by its ability to represent parametric uncertainties in both local and coupling dynamics. This problem class covers a wide range of practical applications, e.g. load frequency control of interconnected power systems \cite{ma2014distributed}, energy efficient building climate control \cite{maasoumy2014handling} or financial optimization \cite{primbs2009stochastic}.

In \cite{farina2016distributed,mark2019distributed}, the authors propose distributed stochastic MPCs (DSMPC) for linear systems with additive uncertainty. Each subsystem optimizes its local input by taking the neighboring state sequences as disturbances to reject. Recursive feasibility is guaranteed via a feasibility-based conditional initialization procedure. These approaches tend to be overly conservative due to their decentralized nature. \\
In \cite{dai2016distributed, dai2016distributed2}, the authors propose a DSMPC for linear systems with parameter uncertainty and bounded additive disturbances. Recursive feasibility is guaranteed by permitting at every time step that only one subsystem optimizes its control sequence, while the other subsystems apply the shifted optimal solution. Due to the sequential update these approaches do not scale with the number of subsystems and quickly become intractable. Moreover, the initialization requires a central warm-start solution.

In this paper, we address the aforementioned downsides of the current state-of-the-art and propose a DSMPC for networks of linear systems subject to chance constraints and multiplicative noise. Compared to \cite{dai2016distributed}, we use dual decomposition to obtain a fully parallelizable DMPC problem that scales with the system dimension and can be solved with the consensus distributed Alternating Direction Method of Multipliers method (ADMM) from \cite{boyd2011distributed}. Furthermore, we provide a fully distributed synthesis method for distributed linear feedback controllers and the distributed terminal ingredients. Thus, the MPC synthesis and the online MPC algorithm both do not rely on a central coordination node. Recursive feasibility of the MPC problem is ensured by adopting two alternative control policies, depending whether the MPC optimization problem is infeasible with the measured state or not. A similar approach is presented in \cite{cannon2009model}, where two control policies are adopted, whether the state of the system lies in a so called invariant set with probability $p$ or not.

The paper is organized as follows: Section \ref{sec:problem_formulation} introduces the distributed system dynamics and chance constraints. Section \ref{sec:distributed_controller} introduces the ingredients of the SMPC problem in a centralized fashion, whereas Section \ref{sec:dist_synthesis} is devoted to a corresponding distributed synthesis procedure. In Section \ref{sec:admm} the distributed SMPC problem and the main result is presented, while in Section \ref{sec:example} a numerical example is carried out. For the sake of readability, the proofs of the results can be found in the appendix.
\subsection*{Notation}
Positive definite and semidefinite matrices are indicated as $A>0$ and $A\geq0$, respectively. Given an event $E_1$ we define the probability occurrence as $\text{Pr}(E_1)$ and the conditional probability given $E_2$ as $\text{Pr}(E_1 | E_2)$. For a random variable $w$, we define the expected value and variance as $\mathbb{E}\{w\}$ and $\text{var}\{w\}$. The conditional expected value and variance of $w$ conditional to a random variable $x$ are denoted as $\mathbb{E}\{w | x\}$ and $\text{var}\{w|x\}$. Two random variables $x,y$ that have the same distribution are equal in distribution, denoted by $x \overset{d}{=} y$. A bar above matrices $\bar{P}$ denotes a lifted matrix into the desired dimension. Local matrices are denoted with a sub index, e.g. $A_{ij}$, whereas global matrices are denoted without any sub index. The weighted 2-norm w.r.t. a positive definite matrix $Q = Q^\top$ is $\Vert x \Vert_Q^2 = x^\top Q x$, whereas $\Vert Q \Vert$ denotes the spectral norm of $Q$, i.e. the largest eigenvalue of $Q$.
\section{Distributed systems and chance constraints}
\label{sec:problem_formulation} 
A distributed system is represented as a graph $\mathcal{G}(\mathcal{N}, \mathcal{E})$ with nodes $\mathcal{N}$ and edges $\mathcal{E}$. Each node $i \in \mathcal{N}$ represents a subsystem with local state $x_i \in \mathbb{R}^{n_i}$ and local input $u_i \in \mathbb{R}^{m_i}$. The neighborhood of subsystem $i$ is given by $\mathcal{N}_i = \{j | (i,j) \in \mathcal{E}) \cup \{i \}$ with the stacked column vector $x_{\mathcal{N}_i} = \text{col}_{j \in \mathcal{N}_i}(x_j) \in \mathbb{R}^{n_{\mathcal{N}_i}}$ and $n_{\mathcal{N}_i} = \sum_{j \in \mathcal{N}_i} n_j$. The global state is given by  $x = \text{col}_{j \in \mathcal{N}}(x_j) \in \mathbb{R}^{n}$ with $n = \sum_{j \in \mathcal{N}} n_j$. To simplify the notation throughout the paper, let $W_i \in \{0,1\}^{\mathcal{N}_i \times n}$ and $T_i \in \{0,1\}^{n_i \times n}$ denote lifting matrices, such that $x_{\mathcal{N}_i} = W_i x$ and $x_i = T_i x$. The set of nodes is defined as $\mathcal{N} = \{1, ..., M\} \subseteq \mathbb{N}$.

\subsection{Problem setup}
The distributed stochastic linear discrete-time system is given by
\begin{eqnarray}
	x_i(k+1)&=& A_{\mathcal{N}_i} x_{\mathcal{N}_i}(k) + B_i u_i(k) \nonumber \\ &+&[C_{\mathcal{N}_i} x_{\mathcal{N}_i}(k) + D_{i} u_i(k)] w_i(k) \quad \forall i \in \mathcal{N},
	\label{eq:system_dynamics}
\end{eqnarray}
where $A_{\mathcal{N}_i} \in \mathbb{R}^{n_i \times n_{\mathcal{N}_i}}$, $B_i \in \mathbb{R}^{n_i \times m_i}$, $C_{\mathcal{N}_i} \in \mathbb{R}^{n_i \times n_{\mathcal{N}_i}}$ and $D_i\in \mathbb{R}^{n_i \times m_i}$. The stochastic uncertainty $w_i \in \mathbb{R}$ is a zero mean white noise with unit variance and unbounded support.
\begin{assum}[Uncorrelated disturbances]
$\mathbb{E}\{w_i(k) w_j(t)\} = 0$ for all $t, k$ and for all $i \neq j$. \label{assum:uncorrelated_dist}
\end{assum}
For each subsystem $i \in \mathcal{N}$ the local states and inputs are subject to the probabilistic constraints
\begin{subequations}
\begin{eqnarray}
&\text{Pr}&(H_{i,r}^x x_i(k) \leq 1) \geq p^x_{i,r} \: \: r = 1, \ldots, n_{i,r}\\
&\text{Pr}&(H_{i,s}^u u_i(k) \leq 1) \geq p^u_{i,s} \: \: s = 1, \ldots, n_{i,s},
\end{eqnarray} \label{eq:chance_constraints}
\end{subequations}
where $H_{i,r}^x \in \mathbb{R}^{n_{i,r} \times n_{i,x}}$ and $H_{i,s}^u \in \mathbb{R}^{n_{i,s} \times n_{i,u}}$, $p_{i,r}^x$ and $p_{i,s}^u$ are the desired probabilities of constraint satisfaction for the $n_{i,r}$ state and $n_{i,s}$ input halfspace constraints. By stacking the local inputs as $u = \text{col}_{i \in \mathcal{N}}(u_i) \in \mathbb{R}^m$, we can write the global system as 
\begin{eqnarray}
	x(k+1) = A x(k) + B u(k) + [C x(k) + D u(k)] w(k).
	\label{eq:global_dynamics_1}
\end{eqnarray}
In this formulation the global matrices $A,C$ are block-sparse and $B,D$ are block-diagonal. We make the following assumption on stabilizability. 
\begin{assum}
	There exists a structured linear feedback controller of the form
	\begin{eqnarray*}
		u \coloneqq K x = \text{col}_{i\in \mathcal{M}} (K_{\mathcal{N}_i} x_{\mathcal{N}_i}), 
	\end{eqnarray*}
	where $K \in \mathbb{R}^{m \times n}$ and $K_{\mathcal{N}_i} \in  \mathbb{R}^{m_i \times n_{\mathcal{N}_i}} \: \forall i \in \mathcal{M}$, such that the global system \eqref{eq:global_dynamics_1} is asymptotically stable in the mean-square sense \cite{el1995state}.
	\label{assumption:feedback}
\end{assum}
\begin{rem}
	Asymptotic stability in the mean-square sense implies that $\mathbb{E}\{x(k)\} \rightarrow 0$ and $\text{var}\{x(k)\} \rightarrow 0$ as $k \rightarrow \infty$. In Section \ref{sec:dlqr} we provide details on how a controller can be found that fulfills Assumption \ref{assumption:feedback}.
\end{rem} 
\section{Stochastic MPC}
\label{sec:distributed_controller}
In the following, we present the ingredients to formulate a tractable stochastic MPC problem. That is, the stochastic dynamics \eqref{eq:system_dynamics} is represented by its mean and covariance, the control policies $u_i(k)$ are restricted to be affine feedback laws and the chance constraints are approximated via Cantelli's inequality. The terminal ingredients are first stated in a centralized fashion, whereas Section \ref{sec:dist_synthesis} is devoted to a corresponding distributed representation.
\subsection{Predictive mean-variance dynamics}
To distinguish between closed-loop and predicted quantities we introduce the notation $x(t|k), u(t|k)$, which denotes a $t$-step ahead prediction of the state and input $x,u$ made at time step $k$. Let $z(t|k) = \mathbb{E}\{ x(t|k) \: | \: x(k) \}$ and consider the distributed error feedback control law
\begin{eqnarray}
	u_i(t|k) = v_i(t|k) + K_{\mathcal{N}_i} (x_{\mathcal{N}_i}(t|k) - z_{\mathcal{N}_i}(t|k)), \forall i \in \mathcal{N}\label{eq:tube_controller}
\end{eqnarray}
where $K_{\mathcal{N}_i}$ is a structured feedback gain according to Assumption \ref{assumption:feedback} and $v_i(\cdot|k)$ a nominal input sequence obtained as a solution of an MPC optimization problem solved at time $k$.
Note that by assumption $w(t|k) \overset{d}{=} w(t+k)$ is a zero-mean white noise, from which we can deduce that the predictive mean $z_i$ evolves according to 
\begin{eqnarray}
	z_i(t+1|k) = A_{\mathcal{N}_i} z_{\mathcal{N}_i}(t|k) + B_i v_i(t|k) \quad \forall i \in \mathcal{N}.  \label{eq:local_nominal_system} 
\end{eqnarray}	
Now we define $\Sigma_i(t|k) = \text{var}\{ x_i(t|k) - z_i(t|k) \: | \: x(k)\}$. Similar to \cite{primbs2009stochastic,farina2016model}, one can show that the covariance dynamics is governed by 
\begin{eqnarray}
 &&\Sigma_i(t+1|k)= \nonumber \\
  &&\big[C_{\mathcal{N}_i} z_{\mathcal{N}_i}(t|k) + D_i v_i(t|k)\big]\big[C_{\mathcal{N}_i} z_{\mathcal{N}_i}(t|k) + D_i v_i(t|k)\big]^\top  \nonumber \\
   && + C_{\mathcal{N}_i,K} {\Sigma}_{\mathcal{N}_i}(t|k) C^{\top}_{\mathcal{N}_i,K} + A_{\mathcal{N}_i,K} {\Sigma}_{\mathcal{N}_i}(t|k) A_{\mathcal{N}_i,K}^\top,
\label{eq:dist_variance_pred} 
\end{eqnarray}
where $\Sigma_{\mathcal{N}_i}(t|k) = \text{var}\{ x_{\mathcal{N}_i}(t|k) - z_{\mathcal{N}_i}(t|k) \: | \: x(k)\}$, $A_{\mathcal{N}_i,K} = A_{\mathcal{N}_i} + B_i K_{\mathcal{N}_i}$ and $C_{\mathcal{N}_i,K} = C_{\mathcal{N}_i} + D_i K_{\mathcal{N}_i}$. 
A problem occurs if we try to implement \eqref{eq:dist_variance_pred} in a distributed setting. That is, the right hand side uses a block dense matrix $\Sigma_{\mathcal{N}_i}$ to update the matrix $\Sigma_i$. In this formulation the off-diagonal blocks of $\Sigma_{\mathcal{N}_i}$ are not updated and lead to a potential underestimation of the actual covariance matrix by neglecting the uncertain coupling dynamics. \\
To solve the aforementioned issue we define a positive semi-definite block-diagonal matrix $\hat{\Sigma}(t|k) \in \mathbb{R}^{n \times n}$ in such a way that
\begin{eqnarray}
\hat{\Sigma}(t|k) = \begin{bmatrix}
	\hat{\Sigma}_1(t|k) & \cdots & 0 \\
	\vdots 		& \ddots & \vdots \\
	0  		      & \ldots & \hat{\Sigma}_M(t|k)
\end{bmatrix} \geq \Sigma(t|k)
\label{eq:block_diagonal_cov}
\end{eqnarray}
and introduce block-diagonal covariance dynamics
\begin{eqnarray}
 &&\hat{\Sigma}_i(t+1|k)= \nonumber \\
  &&\big[{C}_{\mathcal{N}_i} z_{\mathcal{N}_i}(t|k) + {D}_i v_i(t|k)\big]\big[{C}_{\mathcal{N}_i} z_{\mathcal{N}_i}(t|k) + {D}_i v_i(t|k)\big]^\top  \nonumber \\
   && + \tilde{C}_{\mathcal{N}_i,K} \hat{\Sigma}_{\mathcal{N}_i}(t|k) \tilde{C}^{\top}_{\mathcal{N}_i,K} + \tilde{A}_{\mathcal{N}_i,K} \hat{\Sigma}_{\mathcal{N}_i}(t|k) \tilde{A}_{\mathcal{N}_i,K}^\top,
\label{eq:covariance_update_modified} 
\end{eqnarray}
where $\tilde{A}_{\mathcal{N}_i,K} = \sqrt{|\mathcal{N}_i|} A_{\mathcal{N}_i,K}$ and $\tilde{C}_{\mathcal{N}_i,K} = \sqrt{|\mathcal{N}_i|} C_{\mathcal{N}_i,K}$. The block-diagonal neighborhood covariance $\hat{\Sigma}_{\mathcal{N}_i}$ consists of the blocks $\hat{\Sigma}_j(t|k)$ for all $j \in \mathcal{N}_i$. The following claim can be proven with \cite[Lem. 1]{farina2016distributed}.
\begin{claim}
Let $\hat{\Sigma}(t+1|k)$ be the global block-diagonal matrix consistent of $\hat{\Sigma}_i(t+1|k)$ for all i $\in \mathcal{N}$ given by \eqref{eq:covariance_update_modified}. If $\hat{\Sigma}(t|k) \geq {\Sigma}(t|k)$, then it holds that $\hat{\Sigma}(t+1|k) \geq {\Sigma}(t+1|k)$.
\end{claim} 
As shown by \cite{primbs2009stochastic}, the dynamics \eqref{eq:covariance_update_modified} can be represented as a linear matrix inequality (LMI)
\begin{eqnarray}
&&\begin{bmatrix}
 \scalebox{0.92}{$\hat{\Sigma}_i(t+1|k)$} 	& \begin{bmatrix}
  \hspace{1.7em} \star  & \hspace{2.5em} \star &\hspace{1.3em} \star 
\end{bmatrix}  \\
 \begin{bmatrix}
  \scalebox{0.92}{$(\tilde{A}_{\mathcal{N}_i} \hat{\Sigma}_{\mathcal{N}_i}(t|k) + \tilde{B}_i U_{\mathcal{N}_i})^\top$} \\
  \scalebox{0.92}{$(\tilde{C}_{\mathcal{N}_i} \hat{\Sigma}_{\mathcal{N}_i}(t|k) + \tilde{D}_i U_{\mathcal{N}_i})^\top$} \\
   \scalebox{0.92}{$(C_{\mathcal{N}_i} z_{\mathcal{N}_i}(t|k) + D_i v_i(t|k))^\top$}
\end{bmatrix} & \begin{bmatrix}
 & \scalebox{0.92}{$\hat{\Sigma}_{\mathcal{N}_i}(t|k)$} & 0 & 0 \\
 &0 &  \scalebox{0.92}{$\hat{\Sigma}_{\mathcal{N}_i}(t|k)$} & 0 \\
 & 0 & 0 &  \scalebox{0.92}{$I$}
 \end{bmatrix} 
\end{bmatrix} \nonumber \\  &&\geq 0\quad \forall i\in \mathcal{N}, \label{eq:pred_var_lmi}
\end{eqnarray}
where $U_{\mathcal{N}_i} = K_{\mathcal{N}_i} \hat{\Sigma}_{\mathcal{N}_i}(t|k)$. The symbol $\star$ denotes the corresponding transposed quantity. Note that \eqref{eq:local_nominal_system} and \eqref{eq:pred_var_lmi} involve only local variables, i.e. variables that are accessible by subsystem $i$ through information exchange with neighbors $j \in \mathcal{N}_i$, which is amendable to distributed optimization.
\subsection{Chance constraint reformulation}
The individual chance constraints \eqref{eq:chance_constraints} are implemented as probabilistic approximations via Cantelli's inequality. As reported in \cite{farina2016model}, one can show that the chance constraints \eqref{eq:chance_constraints} for the predicted states and inputs $x(t|k), u(t|k)$ are verified for all $t \geq 0$ if we instead impose the following deterministic constraints
\begin{subequations}
\begin{eqnarray}
&& H_{i,r}^x z_i(t|k) \leq 1 -  f(p^x_{i,r})\sqrt{H_{i,r}^{x,\top} \hat{\Sigma}_i(t|k) H_{i,r}^x} \\
&& H_{i,s}^u v_i(t|k) \leq 1 -  f(p^u_{i,s})\sqrt{H_{i,s}^{u,\top} \hat{\Sigma}^u_i(t|k) H_{i,s}^u} 
\end{eqnarray}
\label{eq:nonlinear_constraint_sqrt}
\end{subequations}
for all $i \in \mathcal{N}$, $r = 1, \ldots, n_{i,r}$ and $s = 1, \ldots, n_{i,s}$, where $\hat{\Sigma}^u_i(t|k) = K_{\mathcal{N}_i}\hat{\Sigma}_{\mathcal{N}_i}(t|k) K_{\mathcal{N}_i}^\top$ and $f(p) = \sqrt{p/(1-p)}$.
\begin{rem}
$f(p) = \sqrt{p/(1-p)}$ is a distribution-independent bound on the inverse cumulative density function (quantile function) of $w$. Knowing the exact distribution of the disturbance $w$, one can tighten the bound by replacing $f(p)$ with the exact quantile function.
\end{rem}
In order to make the nonlinear constraints \eqref{eq:nonlinear_constraint_sqrt} applicable to our linear MPC framework, we follow the line of \cite{farina2013probabilistic} and linearize \eqref{eq:nonlinear_constraint_sqrt}, i.e.
\begin{subequations}
\begin{eqnarray}
&&{H_{i,r}^x z_i(t|k) \leq (1 - 0.5\epsilon) - \eta_{i,r}^x H_{i,r}^x \hat{\Sigma}_i(t|k) H_{i,r}^{x,\top}} \label{eq:state_CC}\\
&&{H_{i,s}^u v_i(t|k) \leq (1 - 0.5\epsilon) - \eta_{i,s}^u H_{i,s}^u \hat{\Sigma}^u_i(t|k) H_{i,s}^{u,\top}} \label{eq:constraint_controls}
\end{eqnarray}
\label{eq:constraints_nonlinear}
\end{subequations}
for all $i \in \mathcal{N}$, $r = 1, \ldots, n_{i,r}$ and $s = 1, \ldots, n_{i,s}$, where $\eta_{i,r}^x = {f(p^x_{i,r})^2}/{(2\epsilon)}$, $\eta_{i,s}^u = {f(p^u_{i,s})^2}/{(2\epsilon)}$ and $\epsilon \in (0,1]$ denotes an additional design parameter. The authors of \cite{farina2016model} suggest to take values for $\epsilon$ in the range $[0.3, 0.7]$.
\subsection{Cost function}
\label{sec:cost_fcn}
We minimize the expected quadratic cost function
\begin{eqnarray}
	&& \scalebox{0.85} {$J = \mathbb{E}\bigg \{ \Vert x(N|k) \Vert^2_{P}  +\displaystyle \sum_{t=0}^{N-1} \bigg ( \Vert x(t|k) \Vert^2_{Q} + \Vert u(t|k) \Vert^2_{R} \bigg )  \bigg| x(k) \bigg \}$},
	\label{eq:global:expected_cost}
\end{eqnarray}
where $Q \geq 0$ and $R > 0$ are block-diagonal weighting matrices and $P$ satisfies the following assumption:
\begin{assum}
\label{assum_terminal_cost}
There exists a terminal cost $V_f(x) = \sum_{i \in \mathcal{N}} \Vert x_i \Vert_{P_{i}}^2 = \Vert x \Vert_{P}^2$ with block-diagonal matrix $P > 0$ and a distributed terminal controller $u = K_f x$, such that 
\begin{eqnarray}
	&&(A + B K_f)^\top P (A + B K_f) + (C + D K_f)^\top P (C + D K_f) \nonumber \\
	&&  + Q + K_f^\top R K_f - P \leq 0.
	\label{eq:assum_cost_dec}
\end{eqnarray}
\end{assum}
\begin{rem}
\label{rem:tube_equals_terminal}
The existence of the terminal cost function given in Assumption \ref{assum_terminal_cost} implies that the controller $u = K_f x$ is mean-square stabilizing for the global system \eqref{eq:global_dynamics_1}. This is the same condition we required for the tube controller gain $K$ by Assumption \ref{assumption:feedback}. Therefore, for simplicity we set $K_f = K$ for the remainder of the paper.
\end{rem} 
To evaluate the expected quadratic cost \eqref{eq:global:expected_cost} analytically we resort to standard procedures in linear quadratic stochastic control and separate $J$ into its mean and variance components, such that $J = J_{m} + J_{v}$. From block-diagonality of $Q,R,P$ and $\hat{\Sigma}$ the cost function \eqref{eq:global:expected_cost} is fully separable, i.e. $\hat{J} = \sum_{i \in \mathcal{N}} ( J_{m,i} + J_{v,i} )$ with
\begin{eqnarray*}
	&&\scalebox{0.92} {$J_{m,i} = \Vert z_i(N|k) \Vert^2_{P_{i}} + \displaystyle \sum_{t=0}^{N-1} \bigg( \Vert z_i(t|k) \Vert^2_{Q_i} + \Vert v_i(t|k) \Vert^2_{R_i}$ \bigg) } \\
	&&\scalebox{0.92} {$J_{v,i} =  \text{tr}( P_{i} \hat{\Sigma}_i(N|k)) + \displaystyle \sum_{t=0}^{N-1} \text{tr}( \bar{Q}_i + K_{\mathcal{N}_i}^\top R_i K_{\mathcal{N}_i}) \hat{\Sigma}_{\mathcal{N}_i}(t|k))$}, 
\end{eqnarray*}
where $\bar{Q}_i = W_i T_i^\top Q_i T_i W_i^\top$ is lifted into $\mathbb{R}^{\mathcal{N}_i}$ and $J \leq \hat{J}$, which is a consequence of \eqref{eq:block_diagonal_cov}.
\subsection{Terminal constraints}
\label{sec:global_terminal_constr}
To ensure stability, terminal constraints are typically imposed at the end of the prediction horizon on both the mean $z(N|k)$ and the covariance $\hat{\Sigma}(N|k)$ \cite{farina2016model}, i.e.
\begin{eqnarray}
&& z(N|k) \in \mathbb{Z}_f \label{eq:global_term_constr}\\
&&\hat{\Sigma}(N|k) \leq \hat{\Sigma}_{f} \label{eq:global_covariance_cnstr}. 
\end{eqnarray}
In view of \eqref{eq:global_term_constr}, we utilize the block-diagonal terminal weight $P$ and define an ellipsoidal terminal set $\mathbb{Z}_f = \{z \in \mathbb{R}^n| z^\top P z \leq \alpha \}$, where $\alpha > 0$ is a scaling factor that renders $\mathbb{Z}_f$ as a sublevel set of $V_f(x)$. Thus, $\mathbb{Z}_f$ is positively invariant for the global nominal system $z(t+1|k) = A z(t|k) + B v(t|k)$ under the terminal controller $v = K z$, i.e. $(A+BK) z \in \mathbb{Z}_f \: \forall z \in \mathbb{Z}_f$.
Regarding \eqref{eq:global_covariance_cnstr}, we need to compute a block-diagonal terminal covariance matrix $\hat{\Sigma}_{f}$ that verifies the steady-state condition
\begin{eqnarray}
	&&(A + B K) \hat{\Sigma}_{f} (A + B K)^\top + (C + D K) \hat{\Sigma}_{f} (C + D K)^\top \nonumber \\
	&&  + (C + D K) \Psi (C + D K)^\top  \leq \hat{\Sigma}_{f}, \label{eq:global_steady_state_cov}
\end{eqnarray}
where the matrix $\Psi$ is defined in such a way that  
\begin{eqnarray}
\Psi > \psi I > z z^\top \geq 0 \quad \forall z \in \mathbb{Z}_f \label{eq:global_PSI_constr}
\end{eqnarray}
for some $\psi > 0$. The following must hold for all $z \in \mathbb{Z}_f$
\begin{subequations}
\begin{eqnarray}
&&{H_{i,r}^x z_i \leq (1 - 0.5\epsilon) - \eta_{i,r}^x H_{i,r}^x \hat{\Sigma}_{f,i} H_{i,r}^{x,\top}} \label{eq:global_term_constraint_states}\\
&&{H_{i,s}^u K_{\mathcal{N}_i} z_{\mathcal{N}_i} \leq (1 - 0.5\epsilon) - \eta_{i,s}^u H_{i,s}^u \hat{\Sigma}_{f,i}^u H_{i,s}^{u,\top}} \label{eq:global_term_constraint_input}
\end{eqnarray}
\label{eq:terminal_global_constraints}
\end{subequations}
for all $i \in \mathcal{N}$, $r = 1, \ldots, n_{i,r}$ and $s = 1, \ldots, n_{i,s}$, where $\hat{\Sigma}_{f,i}^u = K_{\mathcal{N}_i} \hat{\Sigma}_{f,\mathcal{N}_i} K_{\mathcal{N}_i}^\top$.
In view of \eqref{eq:global_steady_state_cov} - \eqref{eq:global_PSI_constr} it is always possible to define a sufficiently small set $\mathbb{Z}_f$, such that for all $z \in \mathbb{Z}_f$ the terminal chance constraints \eqref{eq:terminal_global_constraints} are verified. In fact, the smaller $\mathbb{Z}_f$, the smaller $\Psi$ resulting from \eqref{eq:global_PSI_constr} and hence the smaller $\hat{\Sigma}_f$ resulting from \eqref{eq:global_steady_state_cov}.
\subsection{Initial conditions}
\label{sec:initial}
In this work, we use a reset-based initialization scheme where we consider the initial conditions $(z(0|k),\hat{\Sigma}(0|k) )$ as free decision variables to ensure the fundamental property of recursive feasibility \cite{farina2016model}. We define the feedback mode (S1) as $x(0|k) = \mathbb{E}\{x(k) \: | \: x(k)\}$, which is chosen via the condition $( z(0|k), \hat{\Sigma}(0|k) ) = (x(k), 0)$ whenever possible. Since the disturbance $w$ is unbounded, mode (S1) can lead to infeasibility. Therefore, we define the backup mode (S2) as $x(0|k) = \mathbb{E}\{ x(k) \: | \: x(k-1) \}$, which is enforced by the condition $( z(0|k), \hat{\Sigma}(0|k) ) = ( z(1|k-1), \hat{\Sigma}(1|k-1))$ if (S1) is infeasible. Thus, we have the binary constraint
\begin{eqnarray}
	( z(0|k), &&\hat{\Sigma}(0|k) ) \in \{S1, S2\}. \label{eq:global_init_constraint} 
\end{eqnarray}
\subsection{Central MPC optimization problem}
The following MPC optimization problem is solved at every time instant $k\geq 0$
\begin{subequations}
\begin{eqnarray}
	&&\underset{z,v,\hat{\Sigma}}{\text{min}} \quad \sum_{i=1}^M \left( J_{m,i}(z_i(\cdot|k), v_i(\cdot|k)) + J_{v,i}(\hat{\Sigma}_{\mathcal{N}_i}(\cdot|k)) \right) \label{eq:mpc_problem:cost}\\
		&&\text{s.t.} \quad \eqref{eq:local_nominal_system}, \eqref{eq:covariance_update_modified}, \eqref{eq:constraints_nonlinear} \quad \forall t=0,..., N-1 \hspace{0.98em} \forall i \in \mathcal{N} \\
		&& \quad \quad \hspace{0.4em} \eqref{eq:global_init_constraint}, \quad z(N|k) \in \mathbb{Z}_f,  \quad \hat{\Sigma}(N|k) \in \hat{\Sigma}_{f},
\end{eqnarray}
\label{eq:central_MPC}
\end{subequations} 
for all $i \in \mathcal{N}$, $r = 1, \ldots, n_{i,r}$ and $s = 1, \ldots, n_{i,s}$.
The last challenge in solving \eqref{eq:central_MPC} distributedly is the terminal set $\mathbb{Z}_f$. The main requirement for the application of distributed optimization is that $\mathbb{Z}_f$ is decomposable into $M$ subproblems, each of which only involves the variables $z_{\mathcal{N}_i}$. In Section \ref{sec:structured_terminal_set} we propose a structured terminal set as a Cartesian product of local time-varying sets, which satisfy the decomposability property.
\begin{rem}
Recursive feasibility of the MPC optimization problem \eqref{eq:central_MPC} can be guaranteed by suitably characterizing the probabilistic constraints \eqref{eq:constraints_nonlinear} through the binary initialization constraint \eqref{eq:global_init_constraint}. Thus, the resulting MPC is not a state-feedback controller as such, since the control input depends not only on $x(k)$ but also on $z(1|k-1)$. Furthermore, the two initialization strategies (S1) and (S2) render the probability operator ambiguous. In particular, when (S1) is applied at time $k$, we enforce that Pr$(H_{i,r}^x x_i(k + 1) \leq 1 | x(k)) \geq p^x_{i,r}$, while with (S2) we force that Pr$(H_{i,r}^x x_i(k + 1) \leq 1 | x(k - \tau)) \geq p^x_{i,r}$, where $k-\tau$ denotes the last time at which (S1) was feasible. 
In contrast, \cite{cannon2009model} propose a state-feedback approach that verifies the conditional probability constraint Pr$(H_{i,r}^x x_i(k+1) \leq 1 | x(k)) \geq p^x_{i,r}$ for all times $k \geq 0$.
\end{rem}
\section{Distributed Synthesis}
\label{sec:dist_synthesis}
\subsection{Structured terminal cost and distributed controller}
\label{sec:dlqr}
In this subsection, we discuss the distributed synthesis of a distributed controller and structured terminal cost, both of which satisfy Assumption \ref{assum_terminal_cost}. More specifically, we aim to find local quadratic functions 
\begin{eqnarray*}
&& V_{f,i}(x_i) = x_i^\top P_i x_i \hspace{2.6em} \forall i \in \mathcal{N}\\
&&  \gamma_i(x_{\mathcal{N}_i}) = x_{\mathcal{N}_i}^\top \Gamma_{\mathcal{N}_i} x_{\mathcal{N}_i} \quad \forall i \in \mathcal{N},
\end{eqnarray*}
such that the global cost decrease condition \eqref{eq:assum_cost_dec} holds true. 
Similar to \cite{conte2016distributed} we introduce indefinite relaxation functions $\gamma_i(\cdot)$ to allow the local cost $V_{f,i}(x_i)$ to partially increase, as long as the global cost $V_{f}(x)$ always decreases in expectation. These implications translate to the following inequalities
\begin{subequations}
\begin{eqnarray}
	&& \mathbb{E} \{V_{f,i}(x_i^+)\:| \: x_i \} - V_{f,i}(x_i) + l(x_{\mathcal{N}_i}, K_{\mathcal{N}_i} x_{\mathcal{N}_i}) - \gamma_i(x_{\mathcal{N}_i}) \nonumber \\
	&& \quad \quad  \leq 0 \quad \forall i \in \mathcal{N} \label{eq:implication_decrease} \\
				&&\sum_{i=1}^M \gamma_i(x_{\mathcal{N}_i}) \leq 0, \label{eq:implication_gamma}
\end{eqnarray}
\end{subequations}
where $x_i^+ = A_{\mathcal{N}_i,K} x_{\mathcal{N}_i} + C_{\mathcal{N}_i,K} x_{\mathcal{N}_i} w_i$ and the stage cost is given by $l(x_{\mathcal{N}_i}, K_{\mathcal{N}_i} x_{\mathcal{N}_i}) = x_{\mathcal{N}_i}^\top (\bar{Q}_i + K_{\mathcal{N}_i}^\top R_i K_{\mathcal{N}_i}) x_{\mathcal{N}_i}$. The resolution of the expected value in \eqref{eq:implication_decrease} and considering \eqref{eq:implication_gamma} leads to a set nonlinear matrix inequalities in the states $x_{\mathcal{N}_i}$. Since these inequalities have to hold for all $x_{\mathcal{N}_i}$, we obtain
\begin{subequations}
\begin{eqnarray}
&&A_{\mathcal{N}_i,K}^\top P_i A_{\mathcal{N}_i,K} +  C^{\top}_{\mathcal{N}_i,K} P_i C_{\mathcal{N}_i,K} - \bar{P}_i \leq \nonumber \\
                & & \qquad -  (\bar{Q}_i + K_{\mathcal{N}_i}^\top R_i K_{\mathcal{N}_i}) + \Gamma_{\mathcal{N}_i} \quad \forall i \in \mathcal{N} \label{eq:dist_cost_a} \\
&&\sum_{i=1}^M W^\top_i \Gamma_{\mathcal{N}_i} W_i \leq 0, \label{eq:dist_cost_b}
\end{eqnarray}
\label{eq:dist_cost}
\end{subequations}
where $\bar{P}_{i} = W_i T_i^\top P_i T_i W_i^\top$. Condition \eqref{eq:dist_cost_a} is structured by design, i.e., it is fully distributedly solvable, while \eqref{eq:dist_cost_b} connects all subsystems with a system-wide coupling constraint.
\begin{lem}
\label{lem:term_cost}
Conditions \eqref{eq:dist_cost_a}- \eqref{eq:dist_cost_b} are equivalent to the following set of LMIs
\begin{subequations}
\label{eq:term_cost_lmi}
\begin{eqnarray}
&&\begin{bmatrix}
 \bar{E}_i + F_{\mathcal{N}_i} 	& \begin{bmatrix}
   \star & \: \: \star  & \: \star & \: \star 
\end{bmatrix}  \\
 \begin{bmatrix}
  A_{\mathcal{N}_i} E_{\mathcal{N}_i}  + B_i Y_{\mathcal{N}_i} \\
  C_{\mathcal{N}_i} E_{\mathcal{N}_i}  + D_i Y_{\mathcal{N}_i} \\
   \bar{Q}_i^{1/2} E_{\mathcal{N}_i} \\
   R_i^{1/2} Y_{\mathcal{N}_i} 
\end{bmatrix} & \begin{bmatrix}
 & E_i  & 0   & 0 &  0 \\
 & 0    & E_i & 0 &  0 \\
 & 0    & 0   & I &  0 \\
 & 0    & 0   & 0 &  I
 \end{bmatrix}
\end{bmatrix} \nonumber \\
&& \geq 0 \quad \forall i \in \mathcal{N} \label{eq:term_cost_LMI_a} \\
&& \sum_{i=1}^M W^\top_i F_{\mathcal{N}_i} W_i \leq 0, \label{eq:term_cost_LMI_b}
\end{eqnarray}
\end{subequations}
where $E_i = P_i^{-1}$, $\bar{E}_i = W_i T^\top_i P_i^{-1} T_i W_i^\top$, $E_{\mathcal{N}_i} = W_i E W_i^\top$, $F_{\mathcal{N}_i} = E_{\mathcal{N}_i} \Gamma_{\mathcal{N}_i} E_{\mathcal{N}_i}$ and $Y_{\mathcal{N}_i} = K_{\mathcal{N}_i} E_{\mathcal{N}_i}$.
\end{lem}
Since the above lemma is a direct extension of \cite[Lemma 10]{conte2016distributed}, we omit the proof and refer to \cite{conte2016distributed}. A structured controller matrix $K$ that satisfies \eqref{eq:term_cost_lmi} forms a mean-square stabilizing controller $u = K x$ for system \eqref{eq:global_dynamics_1}.
\begin{rem}
\label{rem:block_diagonal_matrix}
 Note that the matrix $F_{\mathcal{N}_i}$ is block-sparse and is the last hindrance for a fully distributed implementation. Therefore, the authors of \cite{conte2016distributed} propose to use block-diagonal upper bounds $F_{\mathcal{N}_i} \leq S_{\mathcal{N}_i}$ together with a neighbor-to-neighbor coupling constraint $\sum_{j \in \mathcal{N}_i} T_j W_j^\top S_{\mathcal{N}_j} W_j T_j^\top \leq 0$ for all $i \in \mathcal{N}$ to replace the system-wide coupling constraint \eqref{eq:term_cost_LMI_b}.
\end{rem}

\subsection{Distributed terminal covariance matrix}
Consider the block-diagonal covariance matrix $\hat{\Sigma}_f$ from \eqref{eq:global_covariance_cnstr} and define the neighborhood matrices $\hat{\Sigma}_{f,\mathcal{N}_i} = W_i \hat{\Sigma}_f W_i^\top$. The terminal covariance condition is obtained by substituting $\hat{\Sigma}_{f,i}, \hat{\Sigma}_{f,\mathcal{N}_i}$ and the terminal controller $v_i = K_{\mathcal{N}_i} z_{\mathcal{N}_i}$ into \eqref{eq:dist_variance_pred}, which yields 
\begin{eqnarray}
&& \hat{\Sigma}_{f,i} = A_{\mathcal{N}_i,K} \hat{\Sigma}_{f,\mathcal{N}_i} A_{\mathcal{N}_i,K}^\top + C_{\mathcal{N}_i,K} \hat{\Sigma}_{f, \mathcal{N}_i} C^\top_{\mathcal{N}_i,K}  \nonumber \\ && \qquad + C_{\mathcal{N}_i,K} \Psi_{\mathcal{N}_i} C^\top_{\mathcal{N}_i,K} \quad \forall i \in \mathcal{N},
\label{eq:dist_variance} 
\end{eqnarray}
where $\Psi_{\mathcal{N}_i} \in \mathbb{R}^{\mathcal{N}_i \times \mathcal{N}_i}$ is an arbitrary block-diagonal state covariance matrix. If we set $\Psi_{\mathcal{N}_i} = \hat{\Sigma}_{f,\mathcal{N}_i}$, the nonlinear matrix inequality \eqref{eq:dist_variance} is LMI representable. Note that $\Psi_{\mathcal{N}_i} = \hat{\Sigma}_{f,\mathcal{N}_i}$ can be done if $\hat{\Sigma}_{f,\mathcal{N}_i} \geq \psi_i I$, while $\psi = \underset{i \in \mathcal{N}}{\text{min}}(\psi_i)$ satisfies \eqref{eq:global_PSI_constr}. Define $U_{\mathcal{N}_{i}} = K_{\mathcal{N}_i} \hat{\Sigma}_{f,\mathcal{N}_{i}}$, then the inequality version of \eqref{eq:dist_variance} can be cast as the following pair of LMI's via Schur complements
\begin{subequations}
\label{eq:dist_LMI}
\begin{eqnarray}
&&\begin{bmatrix}
 \hat{\Sigma}_{f,i} 	& \begin{bmatrix}
  \quad \: \star & \quad  \quad \star \quad \: 
\end{bmatrix}  \\
 \begin{bmatrix}
  (A_{\mathcal{N}_{i}} \hat{\Sigma}_{f,\mathcal{N}_{i}}  + B_i U_{\mathcal{N}_{i}})^\top\\
  (C_{\mathcal{N}_i} \hat{\Sigma}_{f,\mathcal{N}_{i}}  + D_i U_{\mathcal{N}_{i}})^\top 
\end{bmatrix} & \begin{bmatrix}
 & \hat{\Sigma}_{f,\mathcal{N}_{i}} & 0  \\
 & 0   & \frac{1}{2} \hat{\Sigma}_{f,\mathcal{N}_{i}}  \\
 \end{bmatrix}
\end{bmatrix} \nonumber \\  &&\geq 0,  \label{eq:dist_LMI_var_a} \\
	&&\begin{bmatrix}
	    \hat{\Sigma}_{f,\mathcal{N}_i} & {I} \\
		{I} 	& \frac{1}{\psi_i} I
	\end{bmatrix} \geq 0 \label{eq:dist_LMI_var_b}.
\end{eqnarray}
\end{subequations}
\subsection{A unique terminal controller}
Observe that in \eqref{eq:term_cost_lmi} and \eqref{eq:dist_LMI} a structured feedback matrix $K_{\mathcal{N}_i}$ is present. In order to resolve this ambiguity, we have to pose a unique LMI problem that simultaneously solves both LMI's. A trivial idea is to impose an additional uniqueness constraint $U_{\mathcal{N}_{i}} \hat{\Sigma}_{f,\mathcal{N}_{i}}^{-1} = Y_{\mathcal{N}_{i}} E_{\mathcal{N}_{i}}^{-1}$, which, as already stated in \cite{farina2016model}, would render the problem non convex. A simple (but conservative) way of circumventing the convexity issue is to set $\hat{\Sigma}_{f,\mathcal{N}_{i}} = E_{\mathcal{N}_{i}}$ and $U_{\mathcal{N}_{i}} = Y_{\mathcal{N}_i}$.
\begin{prop}
\label{prop:combined_LMI}
Set $\hat{\Sigma}_{f,\mathcal{N}_{i}} = E_{\mathcal{N}_{i}}$, $U_{\mathcal{N}_{i}} = Y_{\mathcal{N}_i} \: \forall i \in \mathcal{N}$. If the following optimization problem admits a feasible solution
\begin{eqnarray*}
	\text{max}  \:&& \sum_{i=1}^M log(det(E_i))    \\
	 s.t. &&\eqref{eq:term_cost_LMI_a}, \eqref{eq:term_cost_LMI_b},\eqref{eq:dist_LMI_var_a},\eqref{eq:dist_LMI_var_b} \quad \forall i \in \mathcal{N}
\end{eqnarray*}
then $P_i$ for all $i \in \mathcal{N}$ are unique and the volume the 1-level set of $V_f(x) = \sum_{i \in \mathcal{N}} \Vert x_i \Vert_{P_{i}} = \Vert x \Vert_{P}$ is maximized.
\end{prop}
\begin{pf}
For $E_{\mathcal{N}_{i}} = \hat{\Sigma}_{f,\mathcal{N}_{i}}, U_{\mathcal{N}_{i}} = Y_{\mathcal{N}_i}, \forall i \in \mathcal{N}$ the LMIs \eqref{eq:term_cost_LMI_a}, \eqref{eq:term_cost_LMI_b}, \eqref{eq:dist_LMI_var_a}, \eqref{eq:dist_LMI_var_b} are convex in $E_{\mathcal{N}_{i}}$ and $Y_{\mathcal{N}_i}$, therefore the minimizer is unique. The objective $\sum_{i=1}^M log(det(E_i))$ is convex and maximizes the volume of the 1-level set of $V_f(x)$ \cite{boyd1994linear}. \qed
\end{pf}
\begin{rem}
\label{rem:terminal_constraint_synt}
Proposition \ref{prop:combined_LMI} yields a separable terminal cost function with weights $P_i$, terminal controllers $K_{\mathcal{N}_i}$, relaxation functions $\Gamma_{\mathcal{N}_i}$ and terminal covariance matrices $\hat{\Sigma}_{f,i} = E_i$ for all $i \in \mathcal{N}$ that satisfy Assumption \ref{assum_terminal_cost}. Infeasibility of the optimization problem in Proposition \ref{prop:combined_LMI} implies that there exist no distributed stabilizing terminal controller for system \eqref{eq:system_dynamics}. In this case we can set $\hat{\Sigma}_{f,i}, P_i, \Gamma_{\mathcal{N}_i}$ to zero for all $i \in \mathcal{N}$ and resort to a zero terminal constraint strategy. Furthermore, in view of Remark \ref{rem:tube_equals_terminal}, we then have to compute a structured stabilizing tube controller $K$ for the error system \eqref{eq:tube_controller}, e.g. via LMIs \eqref{eq:term_cost_lmi}.
\end{rem}
\subsection{Structured terminal sets}
\label{sec:structured_terminal_set}
In the following we propose a structured global terminal set $\hat{\mathbb{Z}}_f$ that replaces the global terminal set $\mathbb{Z}_f$ in the MPC optimization problem \eqref{eq:central_MPC}, such that it is solvable with distributed optimization.
\begin{defn}[Time-varying terminal sets]
\label{def:time_varying}
Let $\mathbb{Z}_f$ be the global terminal set from Section \ref{sec:global_terminal_constr} and define $\alpha$ such that for all $z \in \mathbb{Z}_f$ the constraints \eqref{eq:global_steady_state_cov} - \eqref{eq:terminal_global_constraints} are verified. Define local time-varying terminal sets as
\begin{eqnarray}
\label{eq:time_varying_set}	\mathbb{Z}_{f,i}(\alpha_i(k)) \coloneqq \{z_i \in \mathbb{R}^{n_i}| z^\top_i P_i z_i \leq \alpha_i(k) \} \: \: \forall i \in \mathcal{N}, 
\end{eqnarray}
where $\alpha_i(k)$ is given by the set dynamics
\begin{eqnarray*}
	\alpha_i(k+1) = \alpha_i(k) + z_{\mathcal{N}_i}^\top(k) \Gamma_{\mathcal{N}_i} z_{\mathcal{N}_i}(k) \quad \forall i \in \mathcal{N},
\end{eqnarray*}
with $\sum_{i=1}^M \alpha_i(0) \leq \alpha$ and $\alpha_i(0) \geq 0 \quad \forall i \in \mathcal{N}$. A structured global terminal set for the MPC optimization problem \eqref{eq:central_MPC} is then defined as
\begin{eqnarray}
\label{eq:global_terminal_regions}
\hat{\mathbb{Z}}_f(\alpha_1(k), \ldots, \alpha_M(k)) \coloneqq \prod_{i=1}^M \mathbb{Z}_{f,i}(\alpha_i(k)) \subseteq \mathbb{Z}_f.
\end{eqnarray}
\end{defn}
%\begin{rem}
%With \cite[Lem. 8]{conte2016distributed} it can be proven that the local system state remains within the local time-varying terminal sets \eqref{eq:time_varying_set} under the terminal control law and with \cite[Lem. 9]{conte2016distributed} that the Cartesian product \eqref{eq:global_terminal_regions} is recursively feasible. Note that $\hat{\mathbb{Z}}_f$ is a Cartesian product of local sets $\mathbb{Z}_{f,i}(\alpha_i(k))$, hence, it is applicable to distributed optimization. 
%\end{rem} 
In the following, we compute the scaling factor $\alpha$ (Definition \ref{def:time_varying}), such that for all $z \in \hat{\mathbb{Z}}_f$ the terminal constraints  \eqref{eq:global_PSI_constr} and \eqref{eq:terminal_global_constraints} are satisfied. We propose the following distributed linear program, which is an extension of the optimization problem from \cite[Sec. 4.2]{conte2016distributed}.
\begin{subequations}
\label{eq:terminal_alpha_optimization}
\begin{eqnarray}
	\alpha = &&\underset{\hat{\alpha}}{\text{max}} \quad \hat{\alpha} \\
		&&\text{s.t.} \hspace{0.5em} \Vert P_{i}^{-\frac{1}{2}} H_{i,r}^{x,\top} \Vert^2 \hat{\alpha} \leq (\tilde{h}_{i,r}^x)^2 \nonumber \\
		&& \hspace{2em}  \forall i \in \mathcal{N}, \quad r = 1, \ldots, n_{i,r}  \label{eq:distributed_LP_1}\\
		&& \hspace{1.9em} \Vert P_{\mathcal{N}_i}^{-\frac{1}{2}} K_{\mathcal{N}_i}^\top H_{i,s}^{u,\top} \Vert^2 \hat{\alpha} \leq (\tilde{h}_{i,s}^u)^2 \nonumber \\
		&& \hspace{2em} \forall i \in \mathcal{N}, \quad s = 1, \ldots, n_{i,s}  \label{eq:distributed_LP_2} \\
		&& \quad \quad \Vert P_i^{-1} \Vert \hat{\alpha} \leq \psi_i \quad \forall i \in \mathcal{N}. \label{eq:distributed_LP_3}
\end{eqnarray}
\end{subequations}
where $\tilde{h}_{i,r}^x = (1 - 0.5\epsilon) - \eta_{i,r}^x H_{i,r}^x P_i^{-1} H_{i,r}^{x,\top}$ and $\tilde{h}_{i,s}^u = (1 - 0.5\epsilon) - \eta_{i,s}^u H_{i,s}^u K_{\mathcal{N}_i} P_{\mathcal{N}_i}^{-1} K_{\mathcal{N}_i}^\top H_{i,s}^{u,\top}$ denote the right hand side of \eqref{eq:terminal_global_constraints} with the terminal covariance matrices from Proposition \ref{prop:combined_LMI}, i.e. $\hat{\Sigma}_{f,i} = P_i^{-1}$ and $\hat{\Sigma}_{f,\mathcal{N}_{i}} = P_{\mathcal{N}_{i}}^{-1}$.
\begin{lem}
\label{lem:terminal_set}
Let Assumption \ref{assum_terminal_cost} hold. The solution of optimization problem \eqref{eq:terminal_alpha_optimization} defines the largest feasible level set $\mathbb{Z}_f = \{z \in \mathbb{R}^n| z^\top P z \leq \alpha \}$.
\end{lem}
Once such a global level set $\mathbb{Z}_f$ with size $\alpha$ is given, the local terminal sets from Definition \ref{def:time_varying} can be initialized according to $\sum_{i=1}^M \alpha_i(0) \leq \alpha$.
\section{Distributed Optimization for DSMPC}
\label{sec:admm}
At this point the central MPC problem \eqref{eq:central_MPC} can be written entirely by means of distributed ingredients, that is, we replace the terminal set $\mathbb{Z}_f$ with the structured terminal set from Def. \ref{def:time_varying} and the terminal covariance matrix $\hat{\Sigma}_f$ with the block-diagonal matrix $P^{-1}$. In what follows we present a standard distributed consensus ADMM formulation \cite{boyd2011distributed} that exploits the structure of the cost function and the constraints.
\subsection{ADMM Algorithm}
Let $\xi$ contains all global predictions of the input, state and covariance sequences. Let $y_i$ consists of the state and covariance sequence of the neighboring subsystems as predicted by subsystem $i$, i.e. $z_{\mathcal{N}_i}^i(\cdot|k)$ and $\hat{\Sigma}_{\mathcal{N}_i}^i(\cdot|k)$, and the predicted input $v_i(\cdot|k)$ over the prediction horizon $N$.
In this way, the same variables are contained as independent decision variables in $y_j$ for all $j \in \mathcal{N}_i$ and in $\xi$. To coordinate the local solutions, a consensus constraint $G_i \xi = y_i \quad \forall i \in \mathcal{N}$ is introduced.
The matrices $G_i$ are mapping operators where each row is a unit vector with elements in $\{0, 1\}$. Thus, the communication graph $\mathcal{G}(\mathcal{N}, \mathcal{E})$ is encoded in the matrices $G_i$ for all $i \in \mathcal{N}$ and the vectors $y_i$ can be understood as local copies of those entries in $\xi$, which affect subsystem $i$. The augmented Lagrangian for the consensus constraint can now be written as
\begin{eqnarray*}
	\mathcal{L}_i(y_i, \xi, \lambda_i) &=& J_i(y_i) \nonumber \\
	  &+& \lambda_i^\top( y_i - G_i \xi) + \frac{\rho}{2} \Vert y_i - G_i \xi \Vert_2^2 \quad \forall i \in \mathcal{N},
\end{eqnarray*}
where $\lambda_i$ is a Lagrange multiplier and $\rho > 0$ an augmentation factor. With the augmented Lagrangian it is now possible to decompose the MPC optimization problem \eqref{eq:central_MPC} into $\mathcal{N}$ local optimization problems
\begin{subequations}
	\begin{eqnarray}
		y_i^+ = &&\underset{y_i}{\arg \min} \quad \mathcal{L}_i(y_i, \xi, \lambda_i) \\
		&&\hspace{1em} \text{s.t.} \quad \eqref{eq:local_nominal_system}, \eqref{eq:covariance_update_modified}, \eqref{eq:constraints_nonlinear} , \quad \forall t=0,..., N-1 \\
		&& \quad \quad \hspace{1.46em} z_i(N|k) \in \mathbb{Z}_{f,i}(\alpha_i(k)) \label{eq:mpc:mean_constraint}\\
		&& \quad \quad \hspace{1.46em} \Sigma_i(N|k) \leq \hat{\Sigma}_{f,i} = P_{i}^{-1} \label{eq:mpc:variance_constraint}\\
		&& \quad \quad \hspace{1.46em} (z_i(0|k), \Sigma_i(0|k)) = (z_{i,0}, \Sigma_{i,0})
	\end{eqnarray}
	\label{eq:local_MPC}
\end{subequations}
for all $i \in \mathcal{N}$, $r = 1, \ldots, n_{i,r}$ and $s = 1, \ldots, n_{i,s}$. We introduce the following notation: $y_i^{j+}$ indicates $y_i^+$ predicted by subsystem $j$, and $\xi_i$ denotes those entries in $\xi$ that are affected by subsystem $i$.
\begin{algorithm}[]
\caption{Consensus ADMM}\label{alg}
\begin{algorithmic}[1]
\State For each subsystem $i \in \mathcal{N}$ in parallel:
\State Initialize $\lambda_i = 0$, $\xi_i = 0$ and $(z_{i,0}, \Sigma_{i,0})$ according to (S1) or (S2)
\Repeat
\State Solve MPC Problem \eqref{eq:local_MPC} and obtain $y_i^{+}$
\State Communicate $y_i^+$ to neighbors $j \in \mathcal{N}_i$
\State $\xi_i^+ = \frac{1}{\vert \mathcal{N}_i\vert} \sum_{j\in \mathcal{N}_i} y_i^{j,+}$
\State Communicate $\xi_i^+$ to neighbors $j \in \mathcal{N}_i$
\State $\lambda_i^+ = \lambda_i + \rho (y_i^+ - G_i \xi^+)$ 
\Until{convergence}
\end{algorithmic}
\end{algorithm} 

\begin{rem}
If the cost $J_i(y_i)$ is closed, proper and convex and the unaugmented Lagrangian, that is $\mathcal{L}_i$ without the last term, has a saddle point, then the residuals $G_i \xi - y_i$ converge asymptotically to $0$ \cite{boyd2011distributed}.
\end{rem}
For practical reasons a simple stopping criterion for the ADMM is implemented 
\begin{eqnarray}
\Vert G_i \xi -y_i \Vert_\infty \leq \epsilon_c, \label{eq:stopping_cond}
\end{eqnarray}
which can efficiently be checked between iterates of Algorithm \ref{alg}.
Based on the ADMM algorithm \ref{alg} we are ready to state the online DSMPC algorithm \ref{alg_online_mpc}, which is executed at every time instant $k \geq 0$.
\begin{algorithm}[t!]
\caption{Online DSMPC}\label{alg_online_mpc}
\begin{algorithmic}[1]
\State Measure local states $x_i(k)$ for all $i \in \mathcal{N}$ and share with neighbors
\State Set $(z_{i,0}, \Sigma_{i,0}) = (x_i(k), 0)$ for all $i \in \mathcal{N}$ and solve Problem \eqref{eq:central_MPC} via Alg. \ref{alg}
\If {infeasibility is detected} 
	\State Set $(z_{i,0}, \Sigma_{i,0}) = (z^i_i(1|k-1), \hat{\Sigma}^i_i(1|k-1))$ for all $i \in \mathcal{N}$ and solve Problem \eqref{eq:central_MPC} via Alg. \ref{alg}
\EndIf
	\State Each system $i \in \mathcal{N}$ applies the optimal control input \hspace{3em} $u_i(k) = v_i(0|k) + K_{\mathcal{N}_i} (x_{\mathcal{N}_i}(k) - z^{i*}_{\mathcal{N}_i}(0|k))$
	\State Each system $i \in \mathcal{N}$ updates the local terminal set with $\alpha_i(k+1) = \alpha_i(k) + (z_{\mathcal{N}_i}^{i*})^\top(N|k) \Gamma_{\mathcal{N}_i} z_{\mathcal{N}_i}^{i*}(N|k)$
\State $k \rightarrow k + 1$ and go to step 1
\end{algorithmic}
\end{algorithm}

\begin{thm}
\label{thm:main_result}
If at time $k=0$ Problem \eqref{eq:central_MPC} admits a feasible solution via Alg. \ref{alg_online_mpc}, then it is recursively feasible, $\mathbb{E}\{\vert \vert x(k) \vert \vert^2_Q\} \rightarrow 0$ as $k\rightarrow \infty$ and the chance constraints \eqref{eq:chance_constraints} are satisfied for all times $k\geq 0$.
\end{thm}
\section{Numerical example}
\label{sec:example}
\begin{figure}[htbp]
	\centering
	\includegraphics[width=0.75\linewidth]{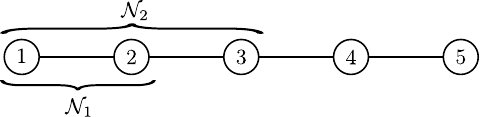}
	\caption{$M=5$ subsystems in a chain topology.} 
	\label{fig:chain}
\end{figure}

In the following section we demonstrate our approach on a numerical example with $M=5$ coupled subsystems in a chain topology, see Figure \ref{fig:chain}.
 For subsystem $1$ it holds $\mathcal{N}_1 = \{1, 2\}$, for subsystem $M$ it holds $\mathcal{N}_M = \{M-1, M\}$ and for subsystems $i = \{2, \ldots, M-1\}$ it holds that $\mathcal{N}_i = \{i-1, i , i+1\}$. 
 Each subsystem $i \in \mathcal{N}$ has the following dynamic matrices
 $A_{ii} = \left[ \begin{smallmatrix}
 	1 & 1\\
 	0 & 1
\end{smallmatrix} \right],
A_{ij} = \left[ \begin{smallmatrix}
	0.1 & 0\\
	0.1 & 0.1
\end{smallmatrix} \right], \forall j \neq i$ and
uncertainty matrices $C_{ii} = \left[ \begin{smallmatrix}
 	0.01 & 0.02\\
 	0.02 & 0.03
\end{smallmatrix} \right], C_{ij} = \left[ \begin{smallmatrix}
	0.002 & 0.02\\
	0    & 0.02
\end{smallmatrix} \right], \forall j \neq i$. The input and uncertain input matrices are given by $
 B_{i} = \left[ \begin{smallmatrix}
  	0\\
 	1
\end{smallmatrix} \right], D_{i} = \left[ \begin{smallmatrix}
	0\\
	0.001
\end{smallmatrix} \right], \forall i \in \mathcal{N}$.
The disturbance is normally distributed with $w \sim \mathcal{N}(0,1)$, the weighting matrices are set to $Q_i = \text{diag}(30,1), R_i = 1, \forall i \in \mathcal{N}$ and the prediction horizon is $N=10$. We consider for each subsystem a single chance constraint $H^x_{1,1} = H^x_{M,1} = [-5 \quad -5]$ and  $H^x_{i,1} = [0 \quad -1.67]$ for all $i \in \{2, \ldots, M-1 \}$, which needs to be satisfied with a probability of at least $p_{i,x} \geq 0.7$. The initial values are set to $x_1(0) = x_M(0) = [3 \quad 0]^\top$ and $x_i(0) = [1 \quad 0]^\top$ for all $i \in \{2, \ldots, M-1\}$. The constraint linearization parameter $\epsilon$ is set to $\epsilon = 0.5$ and $\rho = 10$. \\
In the following, we carry out $N_{mc} = 200$ Monte-carlo simulations to study the closed-loop chance constraint satisfaction and optimality for different parameterizations. In Table \ref{table} we can see for different values of $\epsilon_c$ the average/max number of iterations of Alg. \ref{alg_online_mpc}, the average cost $ av[J] = N_{mc}^{-1} \sum_{q = 1}^{N_{mc}} \sum_{k=0}^{15} \Vert x(k) \Vert_Q^2 + \Vert u(k) \Vert_R^2$, the worst-case empirical in-time constraint satisfaction of subsystem $1$, i.e. $c_{wc} = \!\min_{k \in \{0, \ldots, 15\}} 1 - N_{mc}^{-1}c_v(k)$ with $c_v(k) = \sum_{q = 1}^{N_{mc}} \mathbbm{1}(H^x_{1,1} x_1(k) > 1)$, where $\mathbbm{1}(\cdot)$ denotes the indicator function. Furthermore, we can see the cumulative number of constraint violations for subsystem $1$, that is $c_v = \sum_{q = 1}^{N_{mc}} \sum_{k=0}^{15} c_v(k)$.
\begin{center}
\begin{tabular}{lllllll}
\toprule
$\epsilon_c$ & av{[}it{]} & max{[}it{]} & av{[}J{]} & $c_v$ & $c_{wc}$&\\ \midrule
$5 \cdot 10^{-3}$       & $14.1$     & $38$     & $18627$ & $53$ & 77 \%&\\
$5 \cdot 10^{-4}$       & $23.3$     & $58$     & $18655$ & $57$ & 77 \%&\\
i)        & $-$         & $-$       & $18773$ & $61$ & 76 \%&\\
ii)       & $-$         & $-$       & $18754$ & $49$ & 78 \%&\\
 \bottomrule
\end{tabular} 
\captionof{table}{Impact of $\epsilon_c$ on the performance.}
\label{table}
\end{center}
For comparison we computed two central solutions, where we set up the central SMPC scheme from \cite{farina2016model} according to:
\begin{enumerate}[i)]
	\item the distributed design procedure from this paper 
	\item the centralized design procedure from \cite{farina2016model}. 
\end{enumerate}
If we increase $\epsilon_c$, the average number of iterations increases, which results from the stopping condition \eqref{eq:stopping_cond}. This implies that $\epsilon_c$ directly influences the suboptimality of the solution, i.e., for $\epsilon_c \rightarrow 0$ we restore the optimal solution. Furthermore, it can be seen that for different values $\epsilon_c$, the average cost and the number of cumulative constraint violations vary only slightly compared to the central case i). In each scenario the chance constraints of level $p_{i,x} \geq 0.7$ are empirically verified.
\begin{rem}
The augmentation factor $\rho$ should be selected in appropriate scale to the cost function \eqref{eq:mpc_problem:cost}. If $\rho$ is too small, the primary objective is the minimization of the cost function \eqref{eq:mpc_problem:cost}. As a consequence, the number of iterations until convergence increases. If $\rho$ is too large, the primary objective is the fulfillment of the consensus constraint. Hence, the minimization of the MPC cost becomes less important and the MPC problem \eqref{eq:central_MPC} gets solved inexactly. In practice, finding a good value for $\rho$ is usually done by trial and error.
\end{rem}
\begin{figure}[t!]
\centering
\includegraphics[width=0.9\linewidth]{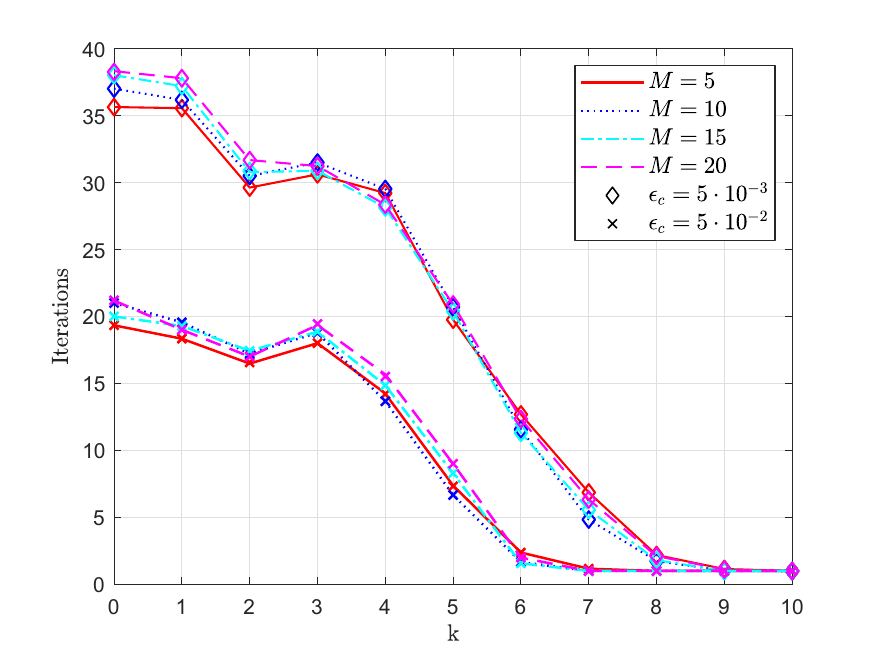}
    \caption{Quantitative impact of the number of subsystems $M$ on the number of iterations of Algorithm \ref{alg}.} 
    \label{fig:quanti}
\end{figure}
In the following, we quantitatively investigate the effect of the number of subsystems $M$ on the online computational demand. In Figure \ref{fig:quanti} we can see for $M \in \{5, 10, 15, 20 \}$ the number of iterations averaged over $100$ Monte-carlo simulations. We consider the same network topology, dynamics, constraints and initial values as before. If we increase the number of subsystems $M$, the number of iterations until convergence slightly increases. This follows from the fact that, the larger we chose the system dimension $M$, the larger the global cost\eqref{eq:mpc_problem:cost} and thus, the more iterations are required to minimize this cost.\\
To summarize: The number of subsystems $M$ affects the number of iterations only marginally, which verifies the scalability of our approach. In order to further reduce the computational demand, we can sacrifice the optimality of the solution by increasing $\epsilon_c$. Since for each test scenario the chance constraints are fulfilled, this boils down to a trade-off between computational complexity and optimality.
\section{Conclusion}
\label{sec:conclusion}
This paper describes a stochastic MPC algorithm for distributed systems with unbounded multiplicative uncertainty. The distributed design guarantees recursive feasibility, point-wise convergence of the states and chance constraint satisfaction. Through the reformulation of the centralized control problem into a distributed semidefinite program, we are able to solve the problem via ADMM. The properties of the controller were highlighted on an example.
\begin{ack}
The authors would like to thank Marcello Farina for his valuable comments on the first draft of the paper.
\end{ack}

{\footnotesize
\bibliographystyle{elsarticle-harv}
\bibliography{automatica_dsmpc_bib}           % and a bib file to produce the 
                                 % bibliography (preferred). The
                                 % correct style is generated by
                                 % Elsevier at the time of printing.
              }                   
\appendix
\section{Proof of Lemma \ref{lem:terminal_set}}
\begin{pf}
Constraints \eqref{eq:distributed_LP_1} and \eqref{eq:distributed_LP_2} are reformulations of \eqref{eq:global_term_constraint_states} and \eqref{eq:global_term_constraint_input} under usage of the support function of the $1$-level set of the elliptical terminal region, see \cite{conte2016distributed} for details. Constraint \eqref{eq:distributed_LP_3} enforces \eqref{eq:global_PSI_constr}, which will be shown in the following. Recall that $\psi = \!\min_{i \in \mathcal{N}}(\psi_i)$ and $P$ is block-diagonal, thus we have the equivalence for all $\alpha \geq 0$
\begin{eqnarray*}
 \Vert P_i^{-1} \Vert {\alpha} \leq \psi_i \quad \forall i \in \mathcal{N} \Longleftrightarrow \Vert P^{-1} \Vert {\alpha} \leq \psi.
\end{eqnarray*}
It remains to show the equivalence of the latter and \eqref{eq:global_PSI_constr}. 
Substitution of $z = P^{-\frac{1}{2}} \tilde{z}$ into the terminal set \eqref{eq:time_varying_set} yields
\begin{eqnarray}
	\forall z \in \mathbb{Z}_{f}: z^\top P z \leq \alpha \Longleftrightarrow \tilde{z}^\top \tilde{z} \leq \alpha \label{eq:alpha_inequality}
\end{eqnarray}
and by substitution into \eqref{eq:global_PSI_constr} that
\begin{eqnarray*}
z z^\top \leq \psi I \Longleftrightarrow P^{-\frac{1}{2}} \tilde{z} \tilde{z}^\top P^{-\frac{1}{2}} \leq \psi I.
\end{eqnarray*}
Taking the norm on both sides yields
\begin{eqnarray*} 
	\Vert P^{-\frac{1}{2}} \tilde{z} \tilde{z}^\top P^{-\frac{1}{2}} \Vert\leq \Vert P^{-1} \Vert \Vert \tilde{z} \tilde{z}^\top \Vert \overset{\eqref{eq:alpha_inequality}}{\leq} \Vert P^{-1} \Vert \alpha \leq \psi,
\end{eqnarray*}
where the second inequality is due to \eqref{eq:alpha_inequality} and the rank one matrix $\tilde{z} \tilde{z}^\top$, which implies that $\Vert \tilde{z} \tilde{z}^\top \Vert = \tilde{z}^\top \tilde{z}$. Since all constraints are convex, maximization of $\hat{\alpha}$ yields the largest feasible level set $\mathbb{Z}_f$. \qed
\end{pf}

\section{Proof of Theorem \ref{thm:main_result}}
%\vspace{-2em}
\begin{pf}
	The proof is inspired by \cite{farina2016model}. Assume that at time $k$ a feasible solution to Problem \eqref{eq:central_MPC} is available. First we prove that at time $k+1$ a feasible solution to Problem \eqref{eq:central_MPC} in mode (S2) exists. To this end, consider the shifted optimal solutions $\scriptstyle \tilde{v}_i(t|k+1) = [v_i^*(1|k), ..., v_i^*(N-1|k), K_{\mathcal{N}_i} z_{\mathcal{N}_i}^*(N|k)]$, $\scriptstyle \tilde{z}_i(t|k+1) = [z_i^*(1|k), ..., z_i^*(N|k), z_i(N+1|k)]$ with $\scriptstyle z_i(N+1|k) = A_{\mathcal{N}_i,K} z^*_{\mathcal{N}_i}(N|k)$ and $\scriptstyle \tilde{\Sigma}_i(t|k+1) = [\hat{\Sigma}_i^*(1|k), ..., \hat{\Sigma}_i^*(N|k), \hat{\Sigma}_i(N+1|k)]$.
From feasibility at time $k$ follows that the state and input constraints \eqref{eq:state_CC}, \eqref{eq:constraint_controls} are verified for any pair $\scriptstyle (\tilde{z}_i(t|k+1), \tilde{\Sigma}_i(t|k+1))$ and $\scriptstyle(\tilde{v}_i(t|k+1), \tilde{\Sigma}_i(t|k+1))$ for each $\scriptstyle t = 0, \ldots, N-1$. At time $\scriptstyle t=N$, we have in view of \eqref{eq:mpc:mean_constraint} and the invariance property of $\scriptstyle \mathbb{Z}_{f,i}(\alpha_i(k))$ (Def. \ref{def:time_varying} and \cite[Lem. 8]{conte2016distributed}) that $\scriptstyle \tilde{z}_i(N|k+1) = z^*_i(N+1|k) \in \mathbb{Z}_{f,i}(\alpha_i)$ and $\scriptstyle \hat{\mathbb{Z}}_f(\alpha_1(k), \ldots, \alpha_M(k)) \subseteq \mathbb{Z}_f$ \cite[Lem. 9]{conte2016distributed}. In view of \eqref{eq:covariance_update_modified}, \eqref{eq:dist_variance} and \eqref{eq:mpc:variance_constraint}, we have that $\scriptstyle \hat{\Sigma}_i(N|k+1) = \hat{\Sigma}_i^*(N+1|k) = \tilde{C}_{\mathcal{N}_i,K} \hat{\Sigma}^*_{\mathcal{N}_i}(N|k) \tilde{C}^{\top}_{\mathcal{N}_i,K} + \tilde{A}_{\mathcal{N}_i,K} \hat{\Sigma}^*_{\mathcal{N}_i}(N|k) \tilde{A}_{\mathcal{N}_i,K}^\top + \tilde{C}_{\mathcal{N}_i,K} z^*_{\mathcal{N}_i}(N|k) (z^*_{\mathcal{N}_i}(N|k))^\top \tilde{C}^{\top}_{\mathcal{N}_i,K} \leq A_{\mathcal{N}_i,K} \hat{\Sigma}_{f,\mathcal{N}_i} A_{\mathcal{N}_i,K}^\top + C_{\mathcal{N}_i,K} \hat{\Sigma}_{f, \mathcal{N}_i} C^\top_{\mathcal{N}_i,K} + C_{\mathcal{N}_i,K} \hat{\Sigma}_{f, \mathcal{N}_i} C^\top_{\mathcal{N}_i,K} \leq \hat{\Sigma}_{f,i}$.
 Hence, both terminal constraints \eqref{eq:mpc:mean_constraint} and \eqref{eq:mpc:variance_constraint} are verified at time $k+1$, which implies satisfaction of the chance constraints \eqref{eq:chance_constraints} for all $k \geq 0$. 

 Next we prove point-wise convergence of the state trajectories. At time step $k+1$ we have to consider the possible shifted optimal solution due to \eqref{eq:global_init_constraint}, i.e. mode (S2). The optimal cost is given by $J^*(k+1) = J_m^*(k+1) + J_v^*(k+1)$ and from optimality follows that $J^*(k+1) \leq J_m(1|k) + J_v(1|k)$ with the suboptimal mean cost
\begingroup
\scriptsize \begin{eqnarray}
&& J_m(1|k) = J_m^*(k) - \sum_{i=1}^M \bigg \{ \vert \vert z_i(0|k) \vert\vert_{Q_i}^2 + \vert\vert v_i^*(0|k) \vert\vert^2_{R_i} \nonumber \\ 
&& - \vert \vert z_i^*(N|k) \vert\vert_{Q_i}^2  - \vert \vert K_{\mathcal{N}_i} z_{\mathcal{N}_i}^*(N|k) \vert\vert_{R_i}^2  \nonumber \\
&& + \vert \vert z_i^*(N|k) \vert\vert_{P_i}^2 - \vert \vert  A_{\mathcal{N}_i,K} z_{\mathcal{N}_i}^*(N|k) \vert\vert_{P_i}^2 \bigg \} \nonumber \\
&& \overset{\eqref{eq:dist_cost_a}}{\leq} J_m^*(k) - \sum_{i=1}^M \bigg \{ \vert \vert z_i(0|k) \vert\vert_{Q_i}^2 \nonumber \\
&& + \vert\vert v_i^*(0|k) \vert\vert^2_{R_i} + \vert \vert z_{\mathcal{N}_i}^*(N|k) \vert \vert^2_{\tilde{W}} - \Vert z_{\mathcal{N}_i}^*(N|k) \Vert^2_{\Gamma_{\mathcal{N}_i}} \bigg \}, \label{eq:mean_cost}
\end{eqnarray}
\endgroup
where $\scriptsize \tilde{W} = C^{\top}_{\mathcal{N}_i,K} P_i C_{\mathcal{N}_i,K}$. Note that $\scriptsize \vert \vert z_i^*(N|k) \vert \vert_{Q_i}^2 = \vert \vert z_{\mathcal{N}_i}^*(N|k) \vert \vert_{\bar{Q}_i}^2$ and $\scriptsize \vert \vert z_i^*(N|k) \vert \vert_{P_i}^2 = \vert \vert z_{\mathcal{N}_i}^*(N|k) \vert \vert_{\bar{P}_i}^2$. The suboptimal variance cost $J_v(1|k)$ is given by
\begingroup
\scriptsize
\begin{eqnarray}
&& J_v(1|k) = J_v^*(k) - \sum_{i=1}^M \bigg\{ \text{tr}(Q_i \hat{\Sigma}_i(0|k)) + \text{tr}(K_{\mathcal{N}_i}^\top R_i K_{\mathcal{N}_i} \hat{\Sigma}_{\mathcal{N}_i}(0|k))  \nonumber \\ 
&& - \text{tr}(Q_i \hat{\Sigma}^*_i(N|k)) - \text{tr}(K_{\mathcal{N}_i}^\top R_i K_{\mathcal{N}_i} \hat{\Sigma}^*_{\mathcal{N}_i}(N|k)) \nonumber \\
&& + \text{tr} \bigg[ P_i \hat{\Sigma}_i^*(N|k) - P_i A_{\mathcal{N}_i,K} \hat{\Sigma}^*_{\mathcal{N}_i}(N|k) A_{\mathcal{N}_i,K}^\top\nonumber \\
&& - P_i (C_{\mathcal{N}_i,K} \hat{\Sigma}_{\mathcal{N}_i}^*(N|k) C_{\mathcal{N}_i,K}^{\top}) \nonumber\\
&&  - P_i (C_{\mathcal{N}_i,K} z_{\mathcal{N}_i}^*(N|k) z_{\mathcal{N}_i}^{*,\top}(N|k) C_{\mathcal{N}_i,K}^{\top})  \bigg] \bigg\} \nonumber \\
&&\overset{\eqref{eq:dist_cost_a}}{\leq} J_v^*(k) - \sum_{i=1}^M \bigg\{ \text{tr}(Q_i \hat{\Sigma}_i(0|k)) + \text{tr}(K_{\mathcal{N}_i}^\top R_i K_{\mathcal{N}_i} \hat{\Sigma}_{\mathcal{N}_i}(0|k) ) \nonumber\\
&& - \vert \vert z^*_{\mathcal{N}_i}(N|k) \vert \vert^2_{\tilde{W}} - \text{tr}(\Gamma_{\mathcal{N}_i} \hat{\Sigma}_{\mathcal{N}_i}^*(N|k)) \big] \bigg\}  \label{eq:var_cost}
\end{eqnarray}
\endgroup
where we used $\scriptstyle \text{tr}(Q_i \hat{\Sigma}^*_i(N|k)) = \text{tr}( \bar{Q}_i \hat{\Sigma}^*_{\mathcal{N}_i}(N|k))$, $\scriptstyle \text{tr}(P_i \hat{\Sigma}^*_i(N|k)) = \text{tr}( \bar{P}_i \hat{\Sigma}^*_{\mathcal{N}_i}(N|k))$ and the cyclic invariance property of the trace to factor out $\scriptstyle \hat{\Sigma}^*_{\mathcal{N}_i}(N|k)$. Furthermore, note that $\vert \vert z^*_{\mathcal{N}_i}(N|k) \vert \vert^2_{\tilde{W}} = \text{tr}(P_i C_{\mathcal{N}_i,K} z_{\mathcal{N}_i}^*(N|k) z_{\mathcal{N}_i}^{*,\top}(N|k) C_{\mathcal{N}_i,K}^{\top})$. After combining \eqref{eq:mean_cost} and \eqref{eq:var_cost} we obtain
\begin{eqnarray*}
&&J^*(k+1) \leq J^*(k) - \sum_{i=1}^M \bigg ( \mathbb{E}\{ \vert \vert x_i(k) \vert \vert_{Q_i}^2 + \vert \vert u_i(k) \vert \vert_{R_i}^2 \} \\
&& - \mathbb{E}\{ \Vert x_{\mathcal{N}_i} \Vert^2_{\Gamma_{\mathcal{N}_i}} \} \bigg ) \overset{\eqref{eq:implication_gamma}}{\leq} J^*(k) - \mathbb{E}\{ \vert \vert x(k) \vert \vert_{Q}^2 \}.
\end{eqnarray*}
Using standard arguments we conclude that $\mathbb{E}\{ \vert \vert x(k) \vert \vert_{{Q}}^2 \} \rightarrow 0$, as $k \rightarrow \infty$. \qed
\end{pf}
\end{document}